\documentclass[12pt,reqno]{amsart}

\usepackage{amssymb}
\usepackage{verbatim}
\usepackage{ifthen}

\newcommand{\showcomments}{no}      

\newcommand{\Sig}{\Sigma}

\newcommand{\del}{\delta}
\newcommand{\sig}{\sigma}

\newcommand{\Z}{{\mathbb Z}}

\newcommand{\cC}{\mathcal C}

\newcommand{\E}{{\sf E}\,}
\renewcommand{\Pr}{{\sf Pr}\,}

\DeclareMathOperator{\Cay}{Cay}
\DeclareMathOperator{\rk}{rk}

\newcommand{\CP}[2][G]{\Cay_{#1}^+({#2})}

\newcommand{\smin}[1][G]{\sig_{\rm min}(#1)}
\newcommand{\smax}[1][G]{\sig_{\rm max}(#1)}
\newcommand{\srnd}[1][G]{\sig_{\rm rnd}(#1)}

\newcommand{\dmin}[1][G]{\del_{\rm min}(#1)}
\newcommand{\dmax}[1][G]{\del_{\rm max}(#1)}
\newcommand{\drnd}[1][G]{\del_{\rm rnd}(#1)}

\renewcommand{\>}{\rangle}

\newcommand{\stm}{\setminus}
\newcommand{\seq}{\subseteq}
\newcommand{\est}{\varnothing}

\newcommand{\longc}{,\ldots,}
\newcommand{\longp}{+\dotsb+}

\theoremstyle{plain}
\newtheorem{claim}{Claim}
\newtheorem{lemma}{Lemma}
\newtheorem{theorem}{Theorem}
\newtheorem{corollary}{Corollary}
\newtheorem{proposition}{Proposition}

\theoremstyle{definition}
\newtheorem{problem}{Problem}

\newcommand{\refm}[1]{~\ref{m:#1}}

\newcommand{\refl}[1]{~\ref{l:#1}}
\newcommand{\reft}[1]{~\ref{t:#1}}
\newcommand{\refc}[1]{~\ref{c:#1}}
\newcommand{\refp}[1]{~\ref{p:#1}}
\newcommand{\refs}[1]{~\ref{s:#1}}
\newcommand{\refb}[1]{~\cite{b:#1}}
\newcommand{\refe}[1]{~\eqref{e:#1}}

\title{Sums and differences \break along Hamiltonian cycles}
\author{Vsevolod F. Lev}
\address{Department of Mathematics,
         The University of Haifa at Oranim,
         Tivon 36006,
         Israel}
\email{seva@math.haifa.ac.il}

\begin{document}
\baselineskip 16pt

\begin{abstract}
Given a finite abelian group $G$, consider the complete graph on the set of
all elements of $G$. Find a Hamiltonian cycle in this graph and for each pair
of consecutive vertices along the cycle compute their sum. What are the
smallest and the largest possible number of sums that can emerge in this way?
What is the expected number of sums if the cycle is chosen randomly? How the
answers change if an orientation is given to the cycle and differences
(instead of sums) are computed? We give complete solutions to some of these
problems and establish reasonably sharp estimates for the rest.
\end{abstract}

\maketitle

\section{Introduction}

For a finite abelian group $G$, by $\cC(G)$ we denote the set of all
Hamiltonian cycles in the complete digraph on the vertex set $G$; thus,
$\cC(G)$ is empty if $G$ is trivial and $|\cC(G)|=(|G|-1)!$ otherwise. Given
a cycle $C\in\cC(G)$, label each edge $(g_1,g_2)\in C$ with the sum $g_1+g_2$
and consider the set $S(C)\seq G$ of all labels along $C$. Now let
\begin{align*}
  \smax &:= \max \{ |S(C)| \colon C \in \cC(G) \}, \\
  \smin &:= \min \{ |S(C)| \colon C \in \cC(G) \}, \\
\intertext{and assuming that $C\in\cC(G)$ is chosen randomly,}
  \srnd &:= \E(|S(C)|).
\end{align*}
Similarly, labeling each (directed) edge $(g_1,g_2)\in C$ with the difference
$g_2-g_1$, consider the set $D(C)\seq G$ of all labels along $C$ and let
\begin{align*}
  \dmax &:= \max \{ |D(C)| \colon C \in \cC(G) \}, \\
  \dmin &:= \min \{ |D(C)| \colon C \in \cC(G) \}, \\
\intertext{and (chosing $C\in\cC(G)$ at random),}
  \drnd &:= \E(|D(C)|).
\end{align*}

In this paper we find the exact values or establish tight bounds for these
six quantities, for all finite abelian groups $G$.

We assume that $\min\est=\max\est=0$ in the definitions above; thus, if $G$
is trivial, then $\smax=\smin=\dmax=\dmin=0$, while $\srnd$ and $\drnd$ are
undefined. Also, if $|G|=2$, then $\smax=\smin=\dmax=\dmin=\srnd=\drnd=1$.

Occasionally, we will consider Hamiltonian cycles on \emph{subsets} of finite
abelian groups, as well as Hamiltonian \emph{paths} on finite abelian groups
or their subsets. The definitions of $S(C)$ and $D(C)$ are carried without
any modification onto the case where $C$ is a Hamiltonian cycle or path on a
finite subset of an additively written group.

In connection with the quantities $\smax$ and $\dmax$ we will be interested
in Hamiltonian cycles and paths such that all sums (differences) of two
consecutive elements along the cycle or path are pairwise distinct; that is,
$|S(C)|=|A|$ (respectively, $|D(C)|=|A|$) for a cycle and $|S(C)|=|A|-1$
(respectively, $|D(C)|=|A|-1$) for a path on the set $A$. We call such cycles
and paths \emph{rainbow-sum} (respectively, \emph{rainbow-difference}) and
use abbreviations like ``RS-cycle'' or ``RD-path''. Under various names such
cycles and paths have been studied by a number of authors; for details and
references see next section and also comments at the end of Sections
\refs{dmax} and \refs{smax}.

Both cycles and paths on the set $A$ will be written as
 $C=(a_1\longc a_{|A|})$, where the components of $C$ list the elements of
$A$; clearly, each Hamiltonian path on $A$ has a unique representation of
this sort, and each cycle has $|A|$ representations.

We close this section with the list of notation, used below in this paper and
not introduced yet:
\begin{align*}
  \<g\>   \;  - &\text{ the subgroup, generated by the group element $g$;} \\
  \Sig(G) \;  - &\text{ the sum of the elements
                                         of the finite abelian group $G$;} \\
  \rk(G)  \;  - &\text{ the rank of the finite abelian group $G$;} \\
  \Z/m\Z  \;  - &\text{ the group of residues modulo
                                                the positive integer $m$;} \\
  \CP{S}  \;  - &\text{ the addition Cayley graph, induced on the finite
                        abelian group $G$} \\
                &\text{ by its subset $S$ (see next section for the
                        definition).}
\end{align*}


\section{Summary of results}\label{s:summary}

We now briefly discuss our principal results; proofs (mostly of combinatorial
nature), comments, and more results are postponed until Sections
\refs{dmin}--\!\!\refs{srnd}.

The smallest possible number of differences along a Hamiltonian cycle can be
determined precisely.
\begin{theorem}\label{t:dmin}
For any finite abelian group $G$ we have $\dmin=\rk(G)$.
\end{theorem}

The situation with the \emph{largest} possible number of differences is
subtler and for some groups there is still room for improvement.
\begin{theorem}\label{t:dmax}
For any finite non-trivial abelian group $G$ we have
  $$ \dmax \le \begin{cases}
                 |G|-1 &\quad \text{if\: $\Sig(G)\neq 0$;} \\
                 |G|-2 &\quad \text{if\: $\Sig(G)=0$.}
               \end{cases} $$
Indeed, if $G$ is not isomorphic to the direct sum of a group of odd order
and a non-cyclic group of order $8$, then equality is attained.
\end{theorem}
Notice, that the condition $\Sig(G)\neq 0$ means that $G$ has exactly one
involution; equivalently, $G$ has exactly one invariant factor of even order.

The proof of Theorem \reft{dmax} uses results of Gordon \refb{g} and Headley
\refb{h} asserting that (i) if $\Sig(G)\neq 0$, then $G$ possesses an
RD-path; (ii) if $\Sig(G)=0$ and Sylow $2$-subgroup of $G$ is not of order
$8$, then the set of non-zero elements of $G$ possesses an RD-cycle. (These
results are based on earlier work of Friedlander, Gordon, and Miller
\refb{fgm}). The question of whether the condition $\Sig(G)=0$ alone, without
any extra assumptions, ensures the existence of an RD-cycle on the set of
non-zero elements of $G$, to our knowledge is open. Answering it in the
affirmative would show that in the estimate of Theorem \reft{dmax} equality
is actually attained for \emph{all} finite non-trivial abelian groups $G$.

For those exceptional groups, not covered by the second assertion of Theorem
\reft{dmax}, a lower bound for $\dmax$ is immediate from our next result.
\begin{theorem}\label{t:drnd}
For any finite non-trivial abelian group $G$ we have
  $$ \drnd=(1-e^{-1})|G|+O(1), $$
with an absolute implicit constant.
\end{theorem}

Observe that the expression $(1-e^{-1})|G|+O(1)$ is not at all surprising in
this context, giving the expected number of %
pairwise distinct elements of $G$, appearing in the sequence of $|G|$
randomly and independently chosen elements.

For a subset $S$ of a finite abelian group $G$, consider the graph with the
vertex set $G$ and the edge set
 $\{(g',g'')\in G\times G\colon g'+g''\in S\}$. We denote this graph by
$\CP{S}$ and call it the \emph{addition Cayley graph, induced on $G$ by $S$}.
Addition Cayley graphs received very little attention in the literature; we
mention the papers \refb{g}, where the clique number of the random addition
Cayley graph is studied, and \refb{cgw}, where Hamiltonicity of addition
Cayley graphs is investigated in the special case that $S$ does not contain
elements of the form $2g$ with $g\in G$. The latter paper is particularly
relevant in our context: in a somewhat unexpected way, it turns out that the
quantity $\smin$ is tightly related to Hamiltonicity of the graphs $\CP{S}$.
Specifically, if $C\in\cC(G)$, then $C$ is a Hamiltonian cycle in
$\CP{S(C)}$; conversely, if $S\seq G$ and $C$ is a Hamiltonian cycle in
$\CP{S}$, then $S(C)\seq S$. (We identify graphs with the digraphs, obtained
by replacing each undirected edge with the pair of corresponding directed
edges. Thus, for instance, if $|G|=2$ and $S$ contains the non-zero element
of $G$, then $\CP{S}$ is considered Hamiltonian.) It follows that $\smin$ is
the minimum size of a subset $S\seq G$ such that $\CP{S}$ is Hamiltonian. We
remark that Hamiltonicity of ``conventional'' Cayley graphs was intensively
studied and in particular, it is well-known that any connected Cayley graph
on a finite abelian group with at least three elements is Hamiltonian; see
\refb{m}. However, apart from the results of \refb{cgw}, nothing seems to be
known on Hamiltonicity of \emph{addition} Cayley graphs. We establish some
properties of the graphs $\CP{S}$ in Section \refs{smin} and as a corollary
determine the value of $\smin$ precisely if $G$ is of even order, and obtain
reasonable estimates if $G$ is of odd order.
\begin{theorem}\label{t:smin}
Let $G$ be a finite non-trivial abelian group. If $|G|$ is even and $G$ is of
type $(m_1\longc m_{\rk(G)})$, then
  $$ \smin = \begin{cases}
               \rk(G)\   &\text{if $m_1=2$}, \\
               \rk(G)+1\ &\text{if $m_1>2$}.
             \end{cases}
  $$
If $|G|$ is odd, then
  $$ \rk(G)+1\le \smin \le 2\rk(G)+1. $$
\end{theorem}
Theorem \reft{smin} shows that $\smin=2$ if $G$ is cyclic of even order
$|G|\ge 4$, and $\smin\in\{2,3\}$ if $G$ is cyclic of odd order. Indeed, we
were able to find $\smin$ for cyclic groups $G$ of odd order, too.
\begin{theorem}\label{t:sminc}
If $G$ is cyclic of order $|G|\ge 3$, then
  $$ \smin = \begin{cases}
                2\ &\text{if $|G|$ is even}, \\
                3\ &\text{if $|G|$ is odd}.
              \end{cases} $$
\end{theorem}

Computations seem to suggest that if $\rk(G)=2$ and $|G|\neq 9$, then
$\smin=3$. One can speculate that, indeed, $\smin=\rk(G)+1$ for all
non-cyclic finite abelian groups $G$ of odd order, with a ``small'' number of
exceptions.

Our next theorem establishes the largest possible number of sums along a
Hamiltonian cycle.
\begin{theorem}\label{t:smax}
For any finite non-trivial abelian group $G$ we have
  $$ \smax = \begin{cases}
               |G|   &\quad \text{if\; $\Sig(G)=0$ and $G$ is not} \\
                     &\quad \text{an elementary abelian $2$-group;} \\
               |G|-1 &\quad \text{if\; $\Sig(G)\neq 0$;} \\
               |G|-2 &\quad \text{if\; $G$ is an elementary abelian} \\
                     &\quad \text{$2$-group and $|G|>2$.}
               \end{cases} $$
\end{theorem}

The proof of Theorem \reft{smax} is based on (i) a theorem due to Beals,
Gallian, Headley, and Jungreis (see \refb{bghj}) claiming that if $G$ is a
finite abelian group with $\Sig(G)=0$, which is not an elementary $2$-group,
then $G$ possesses an RS-cycle; (ii) a construction of an RS-path on every
finite abelian group $G$ with $\Sig(G)\neq 0$.

It is worth mentioning that Theorem \reft{smax} bears relation with Latin
transversals in Cayley tables, as we now explain. Let $G$ be a finite
abelian group. Does the Cayley table of $G$ have a Latin transversal? In
other words, do there exist two permutations $(g_1'\longc g_{|G|}')$ and
$(g_1''\longc g_{|G|}'')$ of the elements of $G$ such that
 $(g_1'+g_1''\longc g_{|G|}'+g_{|G|}'')$ is also a permutation? It is
easily seen that if the answer is positive, then $\Sig(G)=0$; on the other
hand, it was shown in \refb{p} (see also \refb{s} where this was
independently rediscovered) that this condition is also sufficient. Notice
the connection with Snevily's conjecture \refb{sn}, which is that any square
sub-table of the Cayley table of a finite abelian group of odd order
possesses a Latin transversal. Clearly, $G$ has an RS-cycle if and only if
its Cayley table has a Latin transversal of some special sort; namely, one
with $g_1''=g_2',\,g_2''=g_3'\,\longc\, g_n''=g_1'$.

\begin{theorem}\label{t:srnd}
For any finite non-trivial abelian group $G$ we have
  $$ \srnd=(1-e^{-1})|G|+O(1), $$
with an absolute implicit constant.
\end{theorem}

The reader is urged to compare Theorems \reft{drnd} and~\reft{srnd}.

\section{The minimum number of differences: $\dmin$}\label{s:dmin}

\begin{proof}[Proof of Theorem \reft{dmin}]
The case where $G$ is trivial is immediate and we assume for the rest of the
proof that $|G|\ge 2$.

Let $n:=|G|$ and let $C=(g_1\longc g_n)\in\cC(G)$ be a Hamiltonian cycle on
$G$, written in such a way that $g_1=0$. In view of $g_{i+1}-g_i\in D(C)\
(i=1\longc n-1$), every element of $G$ can be represented as a sum of
elements of $D(C)$; thus $D(C)$ generates $G$ and consequently
 $|D(C)|\ge \rk(G)$. It follows that $\dmin\ge\rk(G)$.

To show that $\dmin\le\rk(G)$ we use induction by $\rk(G)$. If $\rk(G)=1$
then $G$ is cyclic and, identifying it with the group $\Z/n\Z$, we consider
the Hamiltonian cycle $C:=(0,1,2\longc n-1)$; clearly, $D(C)=\{1\}$, which
settles the case $\rk(G)=1$. To complete the proof we show that for any
Hamiltonian cycle $C=(h_1\longc h_n)$ on the $n$-element abelian group $H$
and any integer $m\ge 2$, there is a Hamiltonian cycle $C'$ on the group
$H\oplus(\Z/m\Z)$ with $|D(C')|\le|D(C)|+1$. Indeed, it is immediately
verified that one can choose
\begin{align*}
  & C' := (h_1, h_2 \longc h_n,\quad \\
  & \hspace{0.85in} h_n+1, h_1+1\longc h_{n-1}+1, \\
  & \hspace{1.3in}h_{n-1}+2, h_n+2 \longc h_{n-2}+2,  \\
  & \hspace{3in} \vdots \\
  & \hspace{2.2in} h_2+(m-1), h_3+(m-1) \longc h_1+(m-1) ).
\end{align*}
\end{proof}

\section{The maximum number of differences: $\dmax$}\label{s:dmax}

\begin{proof}[Proof of Theorem \reft{dmax}]
Let $C=(g_1\longc g_{|G|})\in\cC(G)$. Since $0\notin D(C)$, we have
$|D(C)|\le |G|-1$. Moreover, if all non-zero elements of $G$ are represented
in $D(C)$, then exactly one of them, say $g$, is represented twice and
therefore
  $$ 0 = (g_2-g_1)\longp(g_{|G|}-g_{|G|-1})+(g_1-g_{|G|}) = \Sig(G) + g; $$
consequently, in this case $\Sig(G)=-g\neq 0$. The first assertion of the
theorem follows.

To prove the second assertion, notice first that if $G$ is not isomorphic to
the direct sum of a group of odd order and a non-cyclic group of order $8$,
then either $\Sig(G)\neq 0$, or Sylow $2$-subgroup of $G$ is not of order $8$.
If $\Sig(G)\neq 0$ then, as shown in \refb{g}, the group $G$ possesses an
RD-path; closing this path (by joining its first and last elements), we get a
Hamiltonian cycle $C\in\cC(G)$ with $|D(C)|\ge |G|-1$. If $\Sig(G)=0$ then
$G$ does not have an RD-path: otherwise, arguing as above we would obtain
$\dmax\ge|G|-1$ which, as we saw, is wrong. It is shown in \refb{h}, however,
that if $\Sig(G)=0$ and Sylow $2$-subgroup of $G$ is not of order $8$, then
the set of non-zero elements of $G$ possesses an RD-cycle. Choosing
arbitrarily two adjacent elements of this cycle and inserting $0$ between
them, we obtain a Hamiltonian cycle $C\in\cC(G)$ with
 $|D(C)|\ge |G|-2$.
\end{proof}

For a survey of results, related to the existence of RD-path and RD-cycles in
finite abelian groups, see \refb{k} or \refb{o}. We notice that the standard
terminology used in \cite{b:fgm,b:h,b:g,b:k,b:o} and a number of other papers
is distinct from that we use here. Specifically, RD-paths are called
\emph{directed terraces}, and those groups possessing an RD-path are called
\emph{sequenceable}; furthermore, RD-cycles on the set of non-zero group
elements are called \emph{directed R-terraces}, and those groups for which
such an RD-cycle exists are called \emph{R-sequenceable}.

\section{The expected number of differences: $\drnd$}

\begin{proof}[Proof of Theorem \reft{drnd}]
Let $n:=|G|$; clearly, $n\ge 3$ can be assumed without loss of generality.
Representing $D(C)$ as a sum of indicator random variables corresponding to
the non-zero elements of $G$, write
\begin{equation}\label{e:d1}
  \drnd = \sum_{g\in G\stm\{0\}} \Pr\{g\in D(C)\}
           = \frac1{(n-1)!} \sum_{g\in G\stm\{0\}}
                      | \{ C \in \cC(G) \colon g \in D(C) \} |.
\end{equation}

Assuming that $g\in G\stm\{0\}$ is fixed, for each $A\seq G$ let $\cC_A(G)$
denote the set of all cycles $C\in\cC(G)$ such that every element $a\in A$ is
followed along the cycle by the element $a+g$. Observe, that if $A$ contains
a coset of the subgroup $\<g\>$, generated by $g$, then $\cC_A(G)$ is empty,
unless $A=\<g\>=G$ (in which case $\cC_A(G)$ consists of one single cycle,
induced by $g$ on $G$).

\begin{claim}\label{m:cCA}
If $A$ does not contain a coset of $\<g\>$, then $|\cC_A(G)|=(n-|A|-1)!$.
\end{claim}

\begin{proof}
For each $c\in G\stm A$ find the non-negative integer $k$ (depending on $c$)
so that $c-(k+1)g\notin A$ and $c-kg\longc c-g\in A$. Consider all chains of
the form $(c-kg\longc c-g,c)$, for all $c\in G\stm A$. These chains partition
$G$, and for $C\in\cC(G)$ we have $C\in\cC_A(G)$ if and only if $C$ is
composed of these chains, following each other in some order. The claim
follows now since the number of chains is $|G\stm A|=n-|A|$.
\end{proof}

Using Claim \refm{cCA} and the inclusion-exclusion principle, we get
\begin{align}
  | \{ C \in \cC(G) \colon g \in D(C) \} |
       &= \Big| \bigcup_{A\seq G\colon |A|=1} \cC_A(G) \Big| \nonumber \\
       &= \sum_{\est\neq A\seq G} (-1)^{|A|+1} |\cC_A(G)| \nonumber \\
       &= \sum_{j=1}^{n-1} (-1)^{j+1} (n-j-1)! \, N_j + (-1)^{n+1} \tau,
                                                                \label{e:d2}
\end{align}
where $N_j$ is the number of those $A\seq G$ with $|A|=j$ such that $A$
contains no coset of $\<g\>$, and $\tau$ equals $1$ if $\<g\>=G$ and equals
$0$ otherwise.

We now claim that if $d$ is the order of $g$ in $G$ (so that $d\mid n$ and
$d\ge 2$), then
\begin{equation}\label{e:d3}
  N_j = \sum_{0\le i\le j/d} (-1)^i \binom{n/d}i \binom{n-id}{j-id}
\end{equation}
holds for each $j\in[1,n-1]$. Indeed, represent $G$ as a union of $n/d$
cosets of $\<g\>$. There are $\binom{n/d}i$ ways to choose $i$ cosets, and
for any choice of $i\le j/d$ cosets there are $\binom{n-id}{j-id}$ ways to
choose $j-id$ elements from the remaining cosets; our claim follows now by
the inclusion-exclusion principle.

Substituting \refe{d2} and \refe{d3} into \refe{d1}, and for $d\mid n$
letting $K_d$ be the number of elements $g\in G$ of order $d$, we get
\begin{align*}
  \drnd &= \frac1{(n-1)!} \sum_{d\mid n,\,d\ge 2} K_d \\
        &\qquad \times \Bigg( \sum_{j=1}^{n-1} (-1)^{j+1} (n-j-1)!
               \sum_{0\le i\le j/d} (-1)^i \binom{n/d}i \binom{n-id}{j-id}
                 \Bigg) + O(1) \\
        &= \frac1{(n-1)!} \sum_{d\mid n,\,d\ge 2} K_d \\
        &\qquad \times \Bigg( \sum_{0\le i\le (n-1)/d} (-1)^i \binom{n/d}i
           (n-id)! \sum_{j=\max\{1,id\}}^{n-1}
             \frac{(-1)^{j+1}}{(j-id)!(n-j)} \Bigg) + O(1) \\
        &= M + R +O(1),
\end{align*}
where $M$ is the part of the expression obtained for $i=0$, and $R$ is the
remaining part (corresponding to positive values of $i$). The former is not
difficult to compute:
\begin{align*}
  M &= n \sum_{d\mid n,\,d\ge 2} K_d
        \sum_{j=1}^{n-1} \frac{(-1)^{j+1}}{j!\,(n-j)} \\
    &= n \sum_{d\mid n,\,d\ge 2} K_d
        \sum_{j=1}^{n-1} \frac{(-1)^{j+1}}{j!}
                 \, \Big( \frac1n + O\Big(\frac{j^2}{n^2}\Big) \Big) \\
    &= \sum_{d\mid n,\,d\ge 2} K_d \Big( 1-e^{-1}+O\Big(\frac1n\Big) \Big) \\
    &= (1-e^{-1})n + O(1).
\end{align*}
To complete the proof it remains to estimate the remainder term $R$. Clearly,
we have
  $$ |R| \le \frac1{(n-1)!} \sum_{d\mid n,\,d\ge 2} K_d
          \sum_{1\le i<n/d} \binom{n/d}i
              (n-id)! \sum_{j=id}^{n-1} \frac1{(j-id)!\,(n-j)} . $$
Consider the internal sum. If $id\le j\le\min\{id+2,n-1\}$ then
  $$ \frac1{(j-id)!\,(n-j)}\le \frac1{n-j} \le \frac3{n-id}, $$
while for $id+3\le j+1\le n-1$ we have
  $$ \frac{(j-id)!\,(n-j)}{(j+1-id)!\,(n-j-1)} = \frac{n-j}{(j+1-id)(n-j-1)}
                                 \le \frac13\frac{n-j}{n-j-1} \le \frac23, $$
and hence
  $$ \sum_{j=id}^{n-1} \frac1{(j-id)!\,(n-j)} = O\Big(\frac1{n-id}\Big). $$
Thus
  $$ |R| \le \frac1{(n-1)!} \sum_{d\mid n,\,d\ge 2} K_d
                                  \sum_{1\le i<n/d} \binom{n/d}i (n-id-1)! $$
and to estimate the sum over $i$ we observe that the summand corresponding to
$i=1$ is $\frac nd(n-d-1)!$, while for $2\le i+1<n/d$ we have
\begin{multline*}
  \frac{\binom{n/d}{i+1}(n-(i+1)d-1)!}{\binom{n/d}i(n-id-1)!}
                = \frac{n/d-i}{i+1} \cdot \frac1{(n-(i+1)d)\dotsb(n-id-1)} \\
  \le \frac1{(i+1)d} \cdot \frac{n-id}{n-id-1}
                                              \le \frac2{(i+1)d} \le \frac12.
\end{multline*}
It follows that
  $$ \sum_{1\le i<n/d} \binom{n/d}i (n-id-1)! \le \frac{2n}d
                                                 \, (n-d-1)! \le n(n-3)! $$
and therefore
  $$ |R| \le \frac{n(n-3)!}{(n-1)!} \sum_{d\mid n,\,d\ge 2} K_d = O(1), $$
completing the proof.
\end{proof}

\section{The minimum number of sums: $\smin$}\label{s:smin}

In this section we establish some general results on Hamiltonicity of
addition Cayley graphs and as a corollary derive Theorems \reft{smin} and
\reft{sminc}. It is worth reminding that we identify undirected graphs with
the corresponding digraphs so that, for instance, the complete graph on two
vertices is treated as Hamiltonian.

The trivial necessary condition for Hamiltonicity is connectedness.
\begin{proposition}\label{p:sminconnect}
Let $S$ be a subset of the finite abelian group $G$. In order for $\CP{S}$ to
be connected it is necessary and sufficient that one of the following
conditions holds:
\begin{itemize}
\item[(i)]  $S$ is not contained in a coset of a proper subgroup of $G$;
\item[(ii)] $S$ is contained in the non-zero coset of an index $2$ subgroup
  of $G$, but not contained in any other coset.
\end{itemize}
\end{proposition}

\begin{proof}
The cases where $G$ is trivial and where $S$ is empty are easy to check.
Assuming that $G$ is non-trivial and $S$ is non-empty, let $H$ be the
smallest subgroup such that $S$ is contained in a coset of $H$; in other
words $H$ is the subgroup, generated by the difference set $S-S$. Observe now
that the component of $0$ in $\CP{S}$ consists of all those elements of $G$,
representable as $s_1-s_2+s_3-\dotsb+(-1)^{k+1}s_k$ with $k\ge 0$ and
$s_1\longc s_k\in S$; that is, this component is the set $H\cup(S+H)$. Thus
$\CP{S}$ is connected if and only if either $H=G$, or $H$ is a subgroup of
index $2$ and $S\seq G\stm H$.
\end{proof}

We remark that for the special case where $S$ does not contain group elements
of the form $2g\ (g\in G)$, the assertion of Proposition \refp{sminconnect}
is equivalent to \cite[Proposition~2.3]{b:cgw}.

In contrast with the ``conventional'' case, connectedness is not sufficient
for Hamiltonicity of an addition Cayley graph. Say, if
$S=\{s_1,s_2\}\seq\Z/n\Z$, where $n\ge 3$ is an integer and $s_1-s_2$ is
co-prime with $n$, then the corresponding graph is connected, but not
Hamiltonian. This follows from the fact that there are elements $g\in G$ with
either $2g=s_1$, or $2g=s_2$; such elements have just one neighbor in the
graph. Moreover, it can be shown that if $n\equiv 3\pmod 4,\; G=\Z/n\Z$, and
$S=\{0,1,3\}\seq G$, then $\CP{S}$ is $2$-connected, but not Hamiltonian.

\begin{corollary}\label{c:sminrk}
Let $S$ be a subset of the finite abelian group $G$ such that $\CP{S}$ is
Hamiltonian. Then $|S|\ge\rk(G)$ and moreover, if $G$ is of type $(m_1\longc
m_{\rk(G)})$ with $m_1>2$, then indeed $|S|\ge\rk(G)+1$.
\end{corollary}

\begin{proof}
Since $\CP{S}$ is Hamiltonian, it is connected, hence $S$ is not contained in
a proper subgroup of $G$ by Proposition \refp{sminconnect}. It follows that
$S$ generates $G$ and therefore $|S|\ge\rk(G)$.

Assume now that $|S|=\rk(G)$. Fix arbitrarily an element $s\in S$ and let $H$
denote the subgroup of $G$, generated by $S-s$. Since $0\in S-s$, we have
$\rk(H)\le|S-s|-1=\rk(G)-1$ whence $H$ is a \emph{proper} subgroup. By
Proposition \refp{sminconnect} and in view of $S\seq s+H$, the index of $H$
in $G$ is $2$, and $s\notin H$. Writing for brevity $r:=\rk(G)$, fix a
generating subset $\{h_1\longc h_{r-1}\}$ of $H$. Since $2s\in H$, there are
integers $u_1\longc u_{r-1}$ such that $2s=u_1h_1\longp u_{r-1}h_{r-1}$. We
now distinguish two cases.

Suppose first that there is an index $i\in[1,r-1]$ such that $u_i$ and the
order of $h_i$ in $H$ are of distinct parity. To simplify the notation
suppose, furthermore, that $i=1$. Find an integer $t$ so that
$(u_1+2t)h_1=h_1$ and set $h_1':=s+th_1$. We have then
  $$ 2h_1'=(u_1+2t)h_1+u_2h_2\longp u_{r-1}h_{r-1} $$
and it follows that $h_1$, and thus also $s$, are contained in the subgroup
of $G$, generated by $h_1',h_2\longc h_{r-1}$. We conclude that this subgroup
is the whole group $G$ and hence $\rk(G)\le r-1$, a contradiction.

We have shown that for each $i\in[1,r-1]$, the order of $h_i$ in $H$ and
$u_i$ are of the same parity. One derives easily that there exists $h\in H$
such that $2s=2h$ and then $G$ is the direct sum of $H$ and the two-element
subgroup, generated by $s-h$. Consequently, $G\cong H\oplus(\Z/2\Z)$ and
since $\rk(G)>\rk(H)$, in the canonical representation of $H$ all direct
summands are of even order. Thus the canonical representation of $G$ is
obtained from that of $H$ by adding $\Z/2\Z$ as a direct summand, meaning
that $m_1=2$.
\end{proof}

For groups of even order the estimate of Corollary \refc{sminrk} is sharp.
\begin{lemma}\label{l:sminupper}
Let $G$ be a finite abelian group of type $(m_1\longc m_{\rk(G)})$. If $|G|$
is even, then there is a subset $S\seq G$ with
  $$ |S| = \begin{cases}
             \rk(G)\   &\text{if $m_1=2$}, \\
             \rk(G)+1\ &\text{if $m_1>2$}
           \end{cases}
  $$
such that $\CP{S}$ is Hamiltonian.
\end{lemma}

\begin{proof}
Let $H<G$ be an index $2$ subgroup with $\rk(H)=\rk(G)-1$ if $m_1=2$, and
$\rk(H)=\rk(G)$ if $m_1>2$. Fix an element $s\in G\stm H$. By Theorem
\reft{dmin}, there is a Hamiltonian cycle $C=(h_1\longc h_{|H|})\in\cC(H)$
such that $|D(C)|=\rk(H)$. Now
 $C':=(h_1,s-h_1,\;h_2,s-h_2,\;\longc h_{|H|},s-h_{|H|})$ is a Hamiltonian
cycle on $G$ satisfying $|S(C')|=|D(C)|+1$ and the assertion follows.
\end{proof}

\begin{lemma}\label{l:sminodd}
Let $S$ be a finite non-trivial abelian group of odd order. Then there is a
subset $S\seq G$ of size $|S|\le 2\rk(G)+1$ such that $\CP{S}$ is
Hamiltonian.
\end{lemma}

\begin{proof}
We write $r:=\rk(G)$ and use induction by $r$. If $r=1$ then $G$ is cyclic
and, identifying it with the group $\Z/(2n+1)\Z$ with a positive integer $n$,
we consider the Hamiltonian cycle
  $$ C := (0,\; 1,2n,\; 2,2n-1\;\longc\; n,n+1) \in \cC(G). $$
One verifies immediately that $|S(C)|=3$, proving the lemma for $r=1$.

Assuming now that $r\ge 2$, write $G=H\oplus F$ where $H$ and $F$ are
subgroups of $G$ such that $\rk(H)=r-1$ and $F$ is cyclic. Write $|H|=m$ and
$|F|=2n+1$ (so that $m\ge 3$ and $n\ge 1$ are integers) and find, using the
induction hypothesis, a Hamiltonian cycle $C_H=(h_1\longc h_m)\in\cC(H)$ such
that $|S(C_H)|\le 2r-1$. Identifying $F$ with the group $\Z/(2n+1)\Z$,
consider the following $n+1$ paths in the complete graph on the vertex set
$G$:
\begin{align*}
  & (h_1\longc h_m), \\
  & (h_1+1,h_2+2n\longc h_m+1,\; h_1+2n,h_2+1\longc h_m+2n), \\
  &     \qquad\qquad\qquad \vdots \\
  & (h_1+n,h_2+(n+1)\longc h_m+n, \\
  &     \hspace{1.85in} h_1+(n+1),h_2+n\longc h_m+(n+1)).
\end{align*}
Straightforward verification shows that the (cyclic) concatenation of these
paths yields a Hamiltonian cycle $C\in\cC(G)$ with
$S(C)=S(C_H)\cup\{h_m+h_1+1,h_m+h_1+(n+1)\}$. Thus $|S(C)|=|S(C_H)|+2\le
2r+1$ and the assertion follows.
\end{proof}

Theorem \reft{smin} is immediate form Corollary \refc{sminrk}, Lemmas
\refl{sminupper} and \refl{sminodd}, and the remark, preceding the statement
of Theorem \reft{smin} in Section \refs{summary}. To prove Theorem
\reft{sminc} we classify those subsets $S$ of the finite abelian group $G$
with $|S|\le 2$ and such that $\CP{S}$ is Hamiltonian.

\begin{lemma}\label{l:sminS2}
Let $G$ be a finite abelian group with $|G|\ge 3$.
\begin{itemize}
\item[(i)]  If $|S|=1$ then $\CP{S}$ is not Hamiltonian;
\item[(ii)] If $|S|=2$ then $\CP{S}$ is Hamiltonian if and only if the
  difference of the two elements of $S$ generates an index $2$ subgroup of
  $G$, and this subgroup is disjoint with $S$.
\end{itemize}
\end{lemma}

\begin{proof}
The first assertion is immediate; to prove the second assertion assume that
$|S|=2$ and write $n=|G|$ and $S=\{s_1,s_2\}$. If $n$ is odd then there is an
element $g\in G$ with $2g=s_1$; the vertex of $\CP{S}$, corresponding to $g$,
has then just one neighbor, whence $\CP{S}$ is not Hamiltonian. Suppose that
$n$ is even. In this case for $\CP{S}$ to be Hamiltonian it is necessary and
sufficient that in the $n$-element sequence
  $$ (0,s_1,s_2-s_1,s_1-(s_2-s_1),2(s_2-s_1)\longc s_1-(n/2-1) (s_2-s_1) ) $$
all elements are pairwise distinct and the last element equals $s_2$. The
latter condition means that the difference $s_2-s_1$ has order $n/2$ in $G$,
and the former condition reduces then to $G=\{0,s_1\}\oplus H$, where $H$ is
the subgroup, generated by $s_2-s_1$. Equivalently, $H$ is a subgroup of
index $2$, to which neither $s_1$ nor $s_2$ belong.
\end{proof}

\begin{corollary}\label{c:ham2}
For any finite abelian group $G$ with $|G|\ge 3$ we have $\smin\ge 2$.
Equality is attained in the last estimate if and only if $G$ has a cyclic
subgroup of index $2$; that is, either $G$ is cyclic of even order, or $G$ is
of type $(2,m)$ with an even $m\ge 2$.
\end{corollary}

To complete our investigation of the quantity $\smin$ we observe that for
$|G|$ even the assertion of Theorem \reft{sminc} follows from Corollary
\refc{ham2}, or alternatively from the combination of Corollary \refc{sminrk}
and Lemma \refl{sminupper}; for $|G|$ odd it follows by combining Corollary
\refc{ham2} and Lemma \refl{sminodd}.

\section{The maximum number of sums: $\smax$}\label{s:smax}

\begin{proof}[Proof of Theorem \reft{smax}]
If $C=(g_1\longc g_{|G|})$ is an RS-cycle over the finite non-trivial abelian
group $G$, then $\Sig(G)=(g_1+g_2)\longp(g_{|G|}+g_1)=2\Sig(G)$, whence
$\Sig(G)=0$. Thus, if $\Sig(G)\neq 0$ then $\smax\le|G|-1$. Moreover, if $G$
is an elementary abelian $2$-group with $|G|>2$, then by Theorem \reft{dmax}
we have $\smax=\dmax\le |G|-2$ and indeed, equality holds provided that
$|G|>2$ and $|G|\neq 8$. Furthermore, if $G$ is elementary abelian of order
$|G|=8$, then a cycle $C\in\cC(G)$ with $|S(C)|=6$ is easy to construct; say,
if $g_1,g_2,g_3$ are three independent elements of $G$, then one can set
  $$ C := (0,\, g_1,\, g_2,\, g_3,\, g_1+g_2,\,
                                      g_1+g_2+g_3,\, g_1+g_3,\, g_2+g_3). $$
This establishes the upper bound for $\smax$ and shows that equality is
attained if $G$ is an elementary abelian $2$-group with $|G|>2$.

If $\Sig(G)=0$ and $G$ is not an elementary abelian $2$-group, then by
\cite[Theorem 6.6]{b:bghj} the group $G$ possesses an RS-cycle, which proves
the assertion in this case.


Finally, suppose that $\Sig(G)\neq 0$ and therefore $G$ is the direct sum of
a (possibly, trivial) group of odd order and a cyclic group of even order.
Identifying the latter group with the corresponding group of residues, we can
then write $G=H\oplus(\Z/n\Z)$, where $H$ is a group of odd order $k:=|H|$,
and $n=2m$ is a positive even integer. If $k=1$ we let $h_1=0$ (the zero
element of the group $H$), if $k\ge 3$ we fix an RS-cycle $(h_1\longc
h_k)\in\cC(H)$; the existence of such an RS-cycle follows from $\Sig(H)=0$.
Consider the Hamiltonian paths $P'$ and $P''$ on $\Z/n\Z$ defined by
\begin{align*}
  P'  &:= (0,m,\; 1,m+1,\; 2,m+2\; \longc\; m-1,2m-1), \\
  P'' &:= (m,0,\; m+1,1,\; m+2,2\; \longc\; 2m-1,m-1).
\end{align*}
We have $S(P')=S(P'')=(\Z/n\Z)\stm\{m-1\}$, so that $P'$ and $P''$ are
RS-paths. For $i=1\longc n$ let $p_i'$ and $p_i''$ denote the $i$\/th
elements of $P'$ and $P''$, respectively. Observing that
$p_n'+p_1''=p_n''+p_1'=m-1$, one verifies easily that
\begin{align*}
  & P := (h_1+p_1',h_1+p_2'\longc h_1+p_n', \\
  & \hspace{0.9in} h_2+p_1'', h_2+p_2''\longc h_2+p_n'', \\
  & \hspace{2.3in} \vdots \\
  & \hspace{1.8in} h_k+p_1',h_k+p_2'\longc h_k+p_n')
\end{align*}
is a Hamiltonian path on $G$ with $S(P)=G\stm\{h_1+h_k+(m-1)\}$, hence an
RS-path on $G$. Closing it (by joining its last and first elements) we obtain
a Hamiltonian cycle $C\in\cC(G)$ with $|S(C)|\ge|S(P)|=|G|-1$, as wanted.
%
%
\end{proof}
We note that RS-cycles were first introduced in \refb{bghj}, where they are
called \emph{harmonious sequences} and those groups with an RS-cycle are
called \emph{harmonious groups}.


\section{The expected number of sums: $\srnd$}\label{s:srnd}

\begin{proof}[Proof of Theorem \reft{srnd}]
Consider the two subgroups of the group $G$ defined by
 $G_0:=\{g\in G\colon 2g=0\}$ and $2G:=\{2g\colon g\in G\}$. Let $n:=|G|$ and
$n_0:=|G_0|$, so that $|2G|=n/n_0$ in view of the isomorphism
 $G/G_0\cong 2G$. Without loss of generality we assume that $n\ge 3$.

Given a Hamiltonian cycle $C\in\cC(G)$, represent $S(C)$ as a sum of
indicator random variables corresponding to the elements of $G$, and write
\begin{equation}\label{e:s1}
  \srnd = \sum_{g\in G} \Pr\{g\in S(C)\}
           = \frac1{(n-1)!} \sum_{g\in G}
                      | \{ C \in \cC(G) \colon g \in S(C) \} |.
\end{equation}
Assuming that $g\in G$ is fixed, for each $A\seq G$ let $\cC_A(G)$ denote the
set of all cycles $C\in\cC(G)$ such that every $a\in A$ is followed along $C$
by $g-a$. Observe, that if $g=a'+a''$ with some $a',a''\in A$, then $a'$ is
to be followed by $a''$ and vice versa along any cycle $C\in\cC_A(G)$;
consequently, $\cC_A(G)$ is empty in this case. On the other hand, if $g$ has
no representations as $g=a'+a''$ with $a',a''\in A$, then it is easy to see
that $|\cC_A(G)|=(n-|A|-1)!$ (compare with Claim \refm{cCA} in the proof of
Theorem \reft{drnd}). Using the inclusion-exclusion principle we get
\begin{align}
  | \{ C \in \cC(G) \colon g \in S(C) \} |
       &= \Big| \bigcup_{A\seq G\colon |A|=1} \cC_A(G) \Big| \nonumber \\
       &= \sum_{\est\neq A\seq G} (-1)^{|A|+1} |\cC_A(G)| \nonumber \\
       &= \sum_{j=1}^{n-1} (-1)^{j+1} (n-j-1)! \, N_j, \label{e:s2}
\end{align}
where $N_j$ is the number of those $A\seq G$ with $|A|=j$ such that $g$
cannot be represented as indicated above.

If $g\notin 2G$ then $G$ is a disjoint union of pairs of the form $(c,g-c)$
with $c\in G$, and in order for $A$ to have the property  that $a'+a''\neq g$
for all $a',a''\in A$, it is necessary and sufficient that $A$ contains at
most one element out of each pair. Thus, $N_j=\binom{n/2}j\, 2^j$ if $j\le
n/2$, and $N_j=0$ if $j>n/2$ in this case. On the other hand, if $g\in 2C$
then there are $n_0$ representations $g=2c$ with $c\in G$; removing such
elements $c$ from $G$, we can split the remaining $n-n_0$ elements into
$(n-n_0)/2$ pairs as above, and now in order for $A$ to have the property in
question it is necessary and sufficient that $A$ contains at most one element
out of each of these pairs (and no non-paired elements). Therefore, in this
case $N_j=\binom{(n-n_0)/2}j\, 2^j$ if $j\le (n-n_0)/2$, and $N_j=0$ if
$j>(n-n_0)/2$. Using these observations, from \refe{s1} and \refe{s2} we
obtain
\begin{multline}\label{e:s3}
  \srnd = \frac1{(n-1)!} \, \Big( 1-\frac1{n_0} \Big) \, n
           \sum_{1\le j\le n/2} (-1)^{j+1} (n-j-1)! \, \binom{n/2}{j} 2^j \\
          + \frac1{(n-1)!} \, \frac{n}{n_0}
                   \sum_{1\le j\le (n-n_0)/2} (-1)^{j+1} (n-j-1)! \,
                                                  \binom{(n-n_0)/2}{j} 2^j,
\end{multline}
where the first summand is to be dropped if $n$ is odd.

Let now $m$ denote one of the numbers $n/2$ and $(n-n_0)/2$ and suppose that
$m$ is an integer. Notice that for any $j\in[1,m]$ we have
  $$ \frac1{(n-1)!}\,(n-j-1)! = \frac1{(n-j)\dotsb(n-1)} > n^{-j} $$
while, on the other hand,
\begin{multline*}
  \frac1{(n-1)!}\,(n-j-1)! = n^{-j}
      \Big(1-\frac1n\Big)^{-1}\dotsb\Big(1-\frac{j}n\Big)^{-1} \\
         < n^{-j} e^{2 \big( \frac1n\longp\frac{j}n \big) }
           \le n^{-j} e^{2j^2/n} < n^{-j} \big( 1 + O(j^2/n) \big).
\end{multline*}
It follows that
\begin{align*}
  \frac1{(n-1)!}\sum_{j=1}^m (-1)^{j+1} (n-j-1)! \binom{m}j 2^j \\
     &\hspace{-1.5in} = - \sum_{j=1}^m \binom{m}j \Bigg( -\frac2n \Bigg)^j
         + O\Bigg( \frac1n\,
             \sum_{j=1}^m \binom{m}j j^2
                                   \Big( \frac2n \Big)^j \Bigg) \\
     &\hspace{-1.5in} = 1 - \Big(1-\frac2n\Big)^m
         + O\Bigg( \frac1n\,
             \sum_{j=1}^m \frac{j^2}{j!}
                                \Big( \frac{2m}n \Big)^j \Bigg) \\
     &\hspace{-1.5in} = 1 - e^{m\ln(1-2/n)}
         + O\Bigg( \frac1n\, \sum_{j=1}^m \frac{j^2}{j!} \Bigg) \\
     &\hspace{-1.5in} = 1 - e^{-2m/n+O(m/n^2)} + O(1/n) \\
     &\hspace{-1.5in} = 1 - e^{-2m/n}(1+O(m/n^2)) + O(1/n) \\
     &\hspace{-1.5in} = 1 - e^{-2m/n} + O(1/n).
\end{align*}
Now if $n$ is odd then $n_0=1$ and \refe{s3} along with the last computation
give
  $$ \srnd = n \big( 1 - e^{-(n-1)/n} + O(1/n) \big)
                                                 = (1-e^{-1})\, n + O(1), $$
as wanted. Similarly, if $n$ is even then
\begin{align*}
  \srnd &= \Big(1-\frac1{n_0}\Big) n \big( 1-e^{-1}+O(1/n) \big)
               + \frac{n}{n_0} \, \big( 1 - e^{-(n-n_0)/n} + O(1/n) \big) \\
        &= (1-e^{-1})\,n
               + \frac{n}{n_0}\,\big( e^{-1} - e^{-(n-n_0)/n} \big) + O(1) \\
        &= (1-e^{-1})\,n + \frac{1-e^{n_0/n}}{n_0/n}\, e^{-1} + O(1) \\
        &= (1-e^{-1})\,n + O(1),
\end{align*}
completing the proof.
\end{proof}

\ifthenelse{\equal{\showcomments}{yes}}{}{}


\begin{thebibliography}{BGHJ91}
\bibitem[BGHJ91]{b:bghj}
  {\sc R.~Beals, J.A.~Gallian, P.~Headley}, and {\sc D.~Jungreis},
  Harmonious groups,
  \emph{J. Comb. Theory, Series A} {\bf 56} (1991), 223--238.
\bibitem[CGW03]{b:cgw}
  {\sc B.~Cheyne, V.~Gupta}, and {\sc C.~Wheeler},
  Hamilton Cycles in Addition Graphs,
  \emph{Rose-Hulman Undergraduate Math. Journal} {\bf 1} (4) (2003)
  (electronic).
\bibitem[FGM78]{b:fgm}
  {\sc R.J.~Friedlander, B.~Gordon}, and {\sc M.D.~Miller},
  On a group sequencing problem of Ringel,
  \emph{Congr. Numer.} {\bf 21} (1978), 307--321.
\bibitem[G61]{b:g}
  {\sc B.~Gordon},
  Sequences in groups with distinct partial products,
  \emph{Pacific J. Math.} {\bf 11} (1961), 1309--1313.
\bibitem[G05]{b:gr}
  {\sc B.~Green},
  Counting sets with small sumset, and the clique number
  of random Cayley graphs,
  \emph{Combinatorica} {\bf 25} (3) (2005), 307--326.
\bibitem[H94]{b:h}
  {\sc P.~Headley},
  R-sequenceability and R*-sequenceability of Abelian $2$-groups,
  \emph{Discrete Math.} {\bf 131} (1--3) (1994), 345--350.
\bibitem[K96]{b:k}
  {\sc A.D.~Keedwell},
  Complete mappings and sequencings of finite groups,
  in \emph{CRC Handbook of Combinatorial Designs}
  (C.J~Colbourn and J.H~Dinitz eds.), CRC Press 1996, 246--253.
\bibitem[M83]{b:m}
  {\sc D.~Maru\v si\v c},
  Hamiltonian circuits in Cayley graphs,
  \emph{Discrete Math.} {\bf 46} (1) (1983), 49--54.
\bibitem[O05]{b:o}
  {\sc M.A.~Ollis},
  On terraces for abelian groups,
  \emph{Discrete Mathematics} {\bf 305} (2005), 250--263.
\bibitem[P47]{b:p}
  {\sc L.J.~Paige},
  A note on finite Abelian groups,
  \emph{Bull. Amer. Math. Soc.} {\bf 53} (1947), 590--593.
\bibitem[P51]{b:p2}
  {\sc \bysame},
  Complete mappings of finite groups,
  \emph{Pacific J. Math.} {\bf 1} (1951), 111--116.
\bibitem[Sn99]{b:sn}
  {\sc H.S.~Snevily},
  Unsolved Problems: The Cayley Addition Table of $\Z_n$,
  \emph{Amer. Math. Monthly} {\bf 106} (6) (1999), 584--585.
\bibitem[Su74]{b:s}
  {\sc J.~Sur\'anyi},
  A certain way of presenting finite abelian groups [Hungarian],
  \emph{Mat. Lapok} {\bf 23} (1972), 257--259 (1974).
\end{thebibliography}
\end{document}